\documentclass[12pt]{amsart}
\usepackage{amscd}      
\usepackage{amssymb}
\usepackage{amsmath, amsthm, graphics}
\usepackage{xypic}      
\LaTeXdiagrams          
\usepackage[all,v2]{xy}
\xyoption{2cell} \UseAllTwocells \xyoption{frame} \CompileMatrices
\allowdisplaybreaks[3]

1

\usepackage{latexsym}
\usepackage{epsfig}
\usepackage{amsfonts}
\usepackage{enumerate}
\usepackage{times}

\newtheorem{prop}{Proposition}

\newtheorem{theorem}{Theorem}

\newtheorem{lemma}[subsection]{Lemma}

\theoremstyle{definition}

\theoremstyle{remark}

\theoremstyle{remark}

\newcommand{\proj}{\mathbb{P}}

\newcommand{\com}{\mathbb{C}}

\newcommand{\barr}{\overline}

\newcommand{\oh}{{\mathcal{O}}}

\def\scup{\mathbin{\text{\scriptsize$\cup$}}}
\def\scap{\mathbin{\text{\scriptsize$\cap$}}}

\newcommand{\DD}{\Delta}

\newcommand{\rarr}{\rightarrow}

\newcommand{\Z}{\mathcal{Z}}

\newcommand{\Q}{{\mathbb{Q}}}

\newcommand{\ev}{{\text{ev}}}

\def\<{\left\langle}
\def\>{\right\rangle}

\newcommand{\gwii}[1]{\left< \hspace{-2pt} \left< \, #1 \,
        \right>  \hspace{-2pt} \right>_{0}}

\newcommand{\gwih}[2]{ \left< \, #2 \, \right>_{#1}}
\newcommand{\gwiig}[1]{\left< \hspace{-2pt} \left< \, #1 \,
    \right> \hspace{-2pt} \right>_{g}}
\newcommand{\gwiih}[2]{\left< \hspace{-2pt} \left< \, #2 \,
    \right> \hspace{-2pt} \right>_{#1}}

\newcommand{\grava}[1]{\tau_{#1}(\gamma_{\alpha})}
\newcommand{\gravua}[1]{\tau_{#1}(\gamma^{\alpha})}
\newcommand{\gravb}[1]{\tau_{#1}(\gamma_{\beta})}

\newcommand{\graval}[1]{\tau_{#1}(\gamma_{\ell})}
\newcommand{\gravual}[1]{\tau_{#1}(\gamma^{\ell})}

\newcommand{\ga}{\gamma_{\alpha}}
\newcommand{\gua}{\gamma^{\alpha}}
\newcommand{\gb}{\gamma_{\beta}}

\newcommand{\vs}{{\mathcal S}}

\newcommand{\vw}{{\mathcal W}}
\newcommand{\vv}{{\mathcal V}}
\newcommand{\vu}{{\mathcal U}}

\begin{document}

\title{New topological recursion relations}

\author{Xiaobo Liu}

\address{Department of Mathematics \\University Of Notre Dame \\ Notre Dame,
IN\\USA}

\author{Rahul Pandharipande}
\address{Department of Mathematics\\Princeton University\\ Princeton, NJ\\USA}

\begin{abstract}
Simple  boundary expressions for the $k^{th}$ power of the
cotangent line class $\psi_1$ on $\overline{M}_{g,1}$ are found
for $k\geq 2g$. The method is by virtual localization
on the moduli space of maps to $\proj^1$.
As a consequence, nontrivial tautological classes in the kernel of the
boundary push-forward map
$$\iota_*:A^*( \overline{M}_{g,2}) \rightarrow A^*(\overline{M}_{g+1})$$
are constructed. The geometry of genus $g+1$ curves then
provides universal equations in genus $g$ Gromov-Witten theory.
As an application, we prove all the Gromov-Witten identities
conjectured recently by K. Liu and H. Xu.
\end{abstract}

\date{\today}

\maketitle

\setcounter{section}{-1}
\section{Introduction}
\subsection{Tautological classes}
Let ${\overline{M}_{g,n}}$ be the moduli space of
stable curves of
genus $g$ with $n$ marked points.
Let $A^*(\overline{M}_{g,n})$ denote the Chow ring
with ${\mathbb Q}$-coefficients.
The system of
tautological rings is defined in \cite{fp3}
to be the set of smallest $\Q$-subalgebras of the Chow rings,
$$R^*(\overline{M}_{g,n}) \subset A^*(\overline{M}_{g,n}),$$
satisfying the
following two properties:
\begin{enumerate}
\item[(i)] The system is closed under push-forward via
all maps forgetting markings:
$$\pi_*: R^*(\overline{M}_{g,n}) \rarr R^*(\overline{M}_{g,n-1}).$$
\item[(ii)] The system is closed under push-forward via
all gluing maps:
$$\iota_*: R^*(\overline{M}_{g_1,n_1\scup\{*\}})
\otimes_{\Q}
R^*(\overline{M}_{g_2,n_2\scup\{\bullet\}}) \rarr
R^*(\overline{M}_{g_1+g_2, n_1+n_2}),$$
$$\iota_*: R^*(\overline{M}_{g, n\scup\{*,\bullet\}}) \rarr
R^*(\overline{M}_{g+1, n}),$$
with attachments along the markings $*$ and $\bullet$.
\end{enumerate}
Natural algebraic constructions typically yield Chow classes
lying in the tautological ring.  For example, the standard $\psi$,
$\kappa$, and $\lambda$ classes in $A^*(\overline{M}_{g,n})$
all lie in the tautological ring.
The tautological rings also possess a rich conjectural structure,
see \cite{fp2} for a detailed discussion.

\pagebreak

The moduli space $\overline{M}_{g,n}$ admits a stratification
by topological type indexed by decorated
graphs. The normalized stratum closures are simply
quotients of products of simpler moduli spaces of
pointed curves. A {\em descendent stratum class} in $R^*(\overline{M}_{g,n})$
 is a push-forward from a stratum $S$ of a monomial in the cotangent line
classes of the special points{\footnote{The special points
correspond to the $n$ markings and the singularities
of curves  parametrized by the stratum.}} of $S$.

A relation in $R^*(\overline{M}_{g,n})$ among descendent
stratum classes yields a
universal genus $g$ equation{\footnote{A genus $g$ equation
is allowed to involve all genera up to $g$.}} in Gromov-Witten theory by the
splitting axiom. For example, the equivalence of boundary
strata in $\overline{M}_{0,4}$ implies the WDVV equation.
Several other relations have since been found \cite{BP,G1,G2,kimliu}.

Let $g\geq 1$.
Boundary expressions for powers $\psi_1^k\in R^*(\overline{M}_{g,1})$
of the cotangent line class are the most basic
{\em topological recursion relations}.
For $k\geq g$, boundary expressions for $\psi_1^k$ have been
proved to exist  \cite{fp3,ionel}. While the arguments
are constructive, the method in practice is very difficult.
The answers for $k=g$ appear, for low $g$, to be rather complicated.
{\footnote{
Boundary relations in codimension $g$ for certain linear
combinations of Hodge classes appear in \cite{AS}.}}

The results of the paper concern simple boundary expressions
for $\psi_1^k$ for $k\geq 2g$.
The relations have two interesting consequences.
The first is the construction of nontrivial classes in the
kernel of the boundary push-forward map
$$\iota_*:A^*( \overline{M}_{g,2}) \rightarrow A^*(\overline{M}_{g+1}).$$
By the splitting axioms of Gromov-Witten theory in genus $g+1$,
we obtain universal equations in genus $g$ from linear combinations
of descendent
stratum classes in the kernel of $\iota_*$.
The possibility for such Gromov-Witten equations
was anticipated earlier in discussions
with Faber, but a nontrivial example was not found.
The existence of such nontrivial equations now opens the door
to new possibilities. Are there equations in Gromov-Witten theory
in genus $g$ obtained by boundary embeddings in even higher genera?
Are there new equations{\footnote{The Gromov-Witten
equations obtained from relations in $R^*(\overline{M}_{0,n})$
are known by Keel's study \cite{keel}. Getzler
has claimed complete knowledge of relations in 
$R^*(\overline{M}_{1,n})$.}}
waiting to be found in genus 0 and 1?

The second consequence of our new topological recursion
relations is a proof of the Gromov-Witten conjectures of K. Liu and
H. Xu \cite{kliu}. The conjectures are universal
relations in Gromov-Witten theory related to high powers of
the cotangent line classes. We prove all the conjectures made there.

\subsection{Topological recursion}
Let $g\geq 1$.
Let $L_1 \rightarrow \overline{M}_{g,1}$
be the cotangent line bundle with fiber $T^*_{p_1}(C)$ at the moduli point
$[C,p_1]\in \overline{M}_{g,1}$.
Let
$$\psi_1= c_1(L_1) \in A^1(\overline{M}_{g,1})$$
be the cotangent line class.
For a genus splitting $g_1+g_2=g$, let
 $$\iota:
\DD_{1,\emptyset}(g_1,g_2)\stackrel{\sim}{=} \overline{M}_{g_1,2}\times
\overline{M}_{g_2,1} \rightarrow \overline{M}_{g,1}$$
denote the boundary divisor parametrizing reducible
curves
 $$C=C_1\scup  C_2$$
satisfying $g(C_i)=g_i$ with a single meeting point,
$$C_1\scap C_2 = p_\star,$$ and
marking
$p_1\in C_1$.
Let
$$\psi_{\star_1}, \psi_{\star_2} \in A^1\big(\DD_{1,\emptyset}(g_1,g_2)\big)$$
denote the cotangent line classes at the 
point $p_\star$. Here, $\psi_{\star_1}$ is
the cotangent line along $C_1$ and $\psi_{\star_2}$ is
the cotangent line along $C_2$.

\begin{theorem} \label{bbt} For $g\geq 1$ and $r\geq 0$,
$$\psi_1^{2g+r} =
\sum_{g_1+g_2=g,\ g_i>0}\
\ \ \sum_{a+b=2g-1+r}
(-1)^a \  \frac{g_2}{g} \cdot
\iota_*\Big(
\psi_{\star_1}^a \psi_{\star_2}^b \scap [\DD_{1,\emptyset}(g_1,g_2)] \Big)$$
in $A^{2g+r}(\overline{M}_{g,1})$.
\end{theorem}

For $r>g-2$,  both sides of the above relation vanish
for dimension reasons.
Theorem \ref{bbt} is nontrivial only if $g\geq 2$ and
 $0\leq r \leq g-2$. On the right side of the relation,
the marking 1 carries {no cotangent line classes}.

Theorem \ref{bbt} and several similar relations are proved in Sections
\ref{xz}-\ref{nnvv}
using the virtual geometry of the moduli space of stable maps
$\overline{M}_{g,n}(\proj^1,1)$. Special
intersections against
the virtual class $[\overline{M}_{g,n}(\proj^1,1)]^{vir}$
of the moduli space,
known to
vanish for geometric reasons, are evaluated via virtual
localization \cite{gp} and pushed-forward to $\overline{M}_{g,n}$
to obtain relations. The technique was first used in \cite{fp1}.

\subsection{Consequences}
Let $g\geq 1$ and $r\geq 1$. Consider the class
\begin{equation}\label{sssw}
\xi_{g,r}=\sum_{a+b=2g+r}(-1)^a
\psi_1^a \psi_2^b \in A^{2g+r}(\overline{M}_{g,2}).
\end{equation}
Let
$\iota: \overline{M}_{g,2} \rightarrow \overline{M}_{g+1}$
be the irreducible boundary map.
As a corollary of the new topological recursion relations, we prove
the following result in Section \ref{xxzz}.

\begin{theorem}\label{vyt} For $g\geq 1$ and $r\geq 1$,
$\ \iota_*(\xi_{g,r}) = 0 \in A^{2g+r+1}(\overline{M}_{g+1}).$
\end{theorem}

For $r$ odd, the push-forward
$\iota_*(\xi_{g,r})$
is easily seen to vanish by the antisymmetry of the sum \eqref{sssw}.
We view the class $\xi_{g,{r}}$ as an uninteresting
element of the
kernel of
$$\iota_*: R^*(\overline{M}_{g,2}) \rightarrow R^*(\overline{M}_{g+1}).$$
The universal Gromov-Witten relation obtained
from $\iota_*(\xi_{g,{r}})=0$ is trivial in the $r$ odd case.

The $r$ even case is much more subtle. Here,
$\xi_{g,r}$ is a remarkable element. For
$r \leq g-2$,
$$\xi_{g,r} \neq 0 \in A^*(\overline{M}_{g,2})$$
since we can compute
\begin{equation*}
\int_{\overline{M}_{g,2}} \xi_{g,r}\cdot 
\psi_2^{g-2-r} \cap [\Delta_{1,2}(1,g-1)] = 
\int_{\overline{M}_{1,2}} \psi_1^2 \cdot
\int_{\overline{M}_{g-1,2}} \psi_2^{3g-4} =
\frac{1}{24}\cdot \frac{1}{24^{g-1}(g-1)!}.
\end{equation*}
The vanishing of $\iota_*(\xi_{g,r})$ is nontrivial --- not
a consequence of any elementary symmetry.
Hence, the associated Gromov-Witten relation is also nontrivial.

\subsection{Gromov-Witten theory}
Let $X$ be a nonsingular projective variety over $\com$ of
dimension $d$. Let $\{ \gamma_\ell \}$ be a basis of $H^*(X,\com)$ with
Poincar\'e dual classes $\{ \gamma^\ell \}$.
The descendent
Gromov-Witten invariants of $X$ are
$$\big\langle \tau_{k_1}(\gamma_{\ell_1}) 
\ldots \tau_{k_n}(\gamma_{\ell_n})
\big\rangle_{g,\beta}^X = \int_{[\overline{M}_{g,n}(X, \beta)]^{vir}}
\psi_1^{k_1}\scup \ev_1^*(\gamma_{\ell_1}) \cdots  \psi_n^{k_n} \scup
\ev_n^*(\gamma_{\ell_n})$$
where $\psi_i$ are the cotangent line classes and
$$\ev_i: \overline{M}_{g,n}(X,\beta) \rightarrow X$$
are the evaluation maps associated to the markings.

Let $\{t^\ell_k \}$ be a set of variables.
Let $F^X_g$  be the generating function of the genus $g$ descendent
invariants,
\begin{equation*}
F^X_g = \sum_{\beta \in H_2(X,\Z)}
q^\beta \sum_{n\ge 0} \frac{1}{n!} \sum_{\substack{\ell_1\dots \ell_n \\
k_1 \dots k_n}} t_{k_n}^{\ell_n} \dots t_{k_1}^{\ell_1}\  \langle
\tau_{k_1}(\gamma_{\ell_1}) \dots \tau_{k_n}(\gamma_{\ell_n})
\rangle_{g,\beta}^X .
\end{equation*}
Double brackets denote differentiation,
$$\big\langle
\big
\langle \tau_{k_1}(\gamma_{\ell_1}) \dots \tau_{k_n}(\gamma_{\ell_n})
\big\rangle
\big\rangle_g^X
= \frac{\partial}{\partial t^{\ell_1}_{k_1}}
\cdots
\frac{\partial}{\partial t^{\ell_n}_{k_n}} \ F^X_g.$$

The Gromov-Witten equation obtained from
Theorem \ref{vyt} is the following result (trivial unless $r$ is even)
conjectured by K. Liu and H. Xu.

\begin{theorem}  \label{gwwq} For $g\geq 0$ and $r\geq 1$,

$$\sum_{a+b=2g+r}\ \sum_{\ell}\ (-1)^a   \big\langle
\big\langle \tau_{a}(\gamma_\ell)
\tau_b(\gamma^\ell) \big\rangle \big\rangle_g^X = 0\ .$$

\end{theorem}
Theorem \ref{gwwq} and several related Gromov-Witten
equations conjectured by Liu-Xu are proved in Section \ref{xl}.
Proofs in case $g\leq 2$ or $r> g-2$ were obtained earlier
in \cite{xliu2}.

\subsection{Acknowledgments}
We thank C. Faber, D. Maulik, and H. Xu for  conversations about
tautological relations and Gromov-Witten theory. 
X.~L. was partially supported
by NSF grant DMS-0505835. R.~P. was partially
supported by NSF grant  DMS-0500187.

\section{Localization relations}
\subsection{$\com^*$-action}
Let $t$
be the generator of the $\com^*$-equivariant ring of a point,
$$A^*_{\com^*}(\bullet) = \com[t].$$
Let $\com^*$ act on $\proj^1$ with tangent weights  $t,-t$ at
the fixed points $0,\infty\in \proj^1$ respectively.
There is an induced $\com^*$-action on the moduli space of
maps $\overline{M}_{g,n}(\proj^1,1)$. A $\com^*$-equivariant
virtual class
$$[\overline{M}_{g,n}(\proj^1,1)]^{vir} \in
A^{\com^*}_{2g+n}(\overline{M}_{g,n}(\proj^1,1))$$
is obtained.
The $\com^*$-equivariant evaluation maps
$$\ev_i: \overline{M}_{g,n}(\proj^1,1) \rightarrow \proj^1$$
determine $\com^*$-equivariant classes
$$\ev_i^*([0]), \ev_i^*([\infty]) \in
A_{\com^*}^{1}\big(\overline{M}_{g,n}(\proj^1,1)\big).$$

Denote the $\com^*$-equivariant universal curve and universal
map by
$$
\pi: U \rightarrow \overline{M}_{g,n}(\proj^1,1), \ \
f: U \rightarrow \proj^1.$$
There is a unique lifting of the $\com^*$-action to
$$\oh_{\proj^1}(-2)\rightarrow \proj^1$$
with fiber weights to be $-t,t$
over the fixed points $0,\infty\in \proj^1$ respectively.
Let $$B=R^1 \pi_* f^*\big(\oh_{\proj^1}(-2)\big) \rightarrow
\overline{M}_{g,n}(\proj^1,1).$$
The sheaf $B$ is $\com^*$-equivariant and locally free
of rank $g+1$. Let
$$c_g(B) \in A^g_{\com^*}\big(\overline{M}_{g,n}(\proj^1,1)\big)$$
be the $g^{th}$ Chern class.

A branch morphism for stable maps to $\proj^1$ has been
defined in \cite{bp},
$$\text{br}: \overline{M}_{g,n}(\proj^1,1) \rightarrow
\text{Sym}^{2g}(\proj^1).$$
The branch morphism is $\com^*$-equivariant.
Let $H_0\subset \text{Sym}^{2g}(\proj^1)$ denote
the hyperplane of $2g$-tuples incident to $0\in \proj^1$. Since
$H_0$ is $\com^*$-invariant,
$$\text{br}^*([H_0]) \in A^1_{\com^*}\big(\overline{M}_{g,n}(\proj^1,1)\big).$$

The total space of $\oh_{\proj^1}(-2) \rightarrow \proj^1$ is
 well-known to be the resolution of the $A_1$ singularity
$\com^2/\mathbb{Z}_2$ with respect to the
action $$ -(z_1,z_2) \mapsto (-z_1,-z_2).$$
A localization approach to the corresponding
(reduced) Gromov-Witten theory along similar lines
is developed in \cite{dm}.

\label{axx}
\subsection{Proof of Theorem \ref{bbt}}
\label{xz}
We obtain a boundary expression for $\psi_1^{2g+r}
\in R^*(\overline{M}_{g,1})$ by localization relations on
$\overline{M}_{g,1}(\proj^1,1)$.
Let
$$I_{g,r} = \ev_1^*([\infty]^{2+r}) \scup c_g(B) \scup \text{br}^*([H_0])
\in A^{g+r+3}_{\com^*} \big(\overline{M}_{g,1}(\proj^1,1)\big).$$
Since the non-equivariant limit of $[\infty]^2$ is 0, the non-equivariant
limit of $I_{g,r}$ is also 0.
Let
$$\epsilon: \overline{M}_{g,1}(\proj^1,1) \rightarrow \overline{M}_{g,1}$$
be the forgetful map.
The map $\epsilon$ is $\com^*$-equivariant with respect to the
trivial $\com^*$-action on $\overline{M}_{g,1}$.
After push-forward,
\begin{equation}\label{iiir}
\epsilon_*\big(I_{g,r} \scap [\overline{M}_{g,1}(\proj^1,1)]^{vir}\big) \in
A_{\com^*}^{2g+r}(\overline{M}_{g,1}).
\end{equation}
The virtual localization formula \cite{gp} gives an
explicit calculation of \eqref{iiir} in term of tautological classes.
Setting the non-equivariant limit to 0,
\begin{equation}\label{vvp}
\epsilon_*\big(I_{g,r} \scap [\overline{M}_{g,1}
(\proj^1,1)]^{vir}\big)|_{t=0} = 0,
\end{equation}
yields an equation in $R^{2g+r}(\overline{M}_{g,1})$.

The localization computation of \eqref{vvp} is a sum over residue
contributions of the
$\com^*$-fixed loci of $\overline{M}_{g,1}(\proj^1,1)$.
The contributing $\com^*$-fixed loci $\overline{M}^{\com^*}_{g_1,g_2}$
are indexed by
genus splittings $g_1+g_2=g$.
If $g_1,g_2>0$, the $\com^*$-fixed locus is
\begin{equation}
\label{bpw}
\overline{M}^{\com^*}_{g_1,g_2} \stackrel{\sim}{=}
\overline{M}_{g_1,2} \times \overline{M}_{g_2,1}\subset
\overline{M}_{g,1}(\proj^1,1),
\end{equation}
parametrizing maps with collapsed components of genus $g_1,g_2$ over
$\infty,0 \in \proj^1$ respectively and the marking over $\infty$.
The restriction of $\epsilon$ to the locus \eqref{bpw} is isomorphic to
$$\iota: \Delta_{1,\emptyset}(g_1,g_2) \rightarrow \overline{M}_{g,1}.$$
In the degenerate cases
$$(g_1,g_2)=(0,g) \ \ \text{or} \ \ (g,0),$$
the $\com^*$-fixed loci are isomorphic to $\overline{M}_{g,1}$ and
$\overline{M}_{g,2}$ respectively.

By the virtual localization formula, we obtain
$$
\epsilon_*(I_{g,r} \scap [\overline{M}_{g,1}(\proj^1,1)]^{vir}) =
\sum_{g_1+g_2=g, \ g_i \geq0} \ \epsilon_*\Big(
\frac{I_{g,r}}
{e(\text{Norm}_{g_1,g_2}^{vir})} \scap [\overline{M}^{\com^*}_{g_1,g_2}]
\Big).$$
If $g_1,g_2>0$,
the restriction of $B$ to $\Delta_{1,\emptyset}(g_1,g_2)$
is
$$\mathbb{E}^\vee_{g_1} \otimes(+t) \oplus \mathbb{E}^\vee_{g_2} \otimes (-t) \oplus \com$$
where $\mathbb{E}$ denote the Hodge bundle.
The class $\text{br}^*(H_0)$  restricts to $2g_2t$.
The Euler class of the virtual normal bundle is
$$\frac{1}
{e(\text{Norm}^{vir})} =
\frac{c_{g_2}(\mathbb{E}^\vee\otimes(+t)) c_{g_1}(\mathbb{E}^\vee\otimes
(-t))}{-t^2(t-\psi_{\star_2})(-t-\psi_{\star_1})}.$$
Putting all the terms together and using Mumford's relation{\footnote{Mumford's
relation here is $c_g(\mathbb{E}_g^\vee \otimes(+t)) \cdot
c_g(\mathbb{E}_g^\vee \otimes(-t)) = t^{g}(-t)^g$}} twice, we obtain
$$
\epsilon_*\Big(
\frac{I_{g,r}}
{e(\text{Norm}_{g_1,g_2}^{vir})} \scap [\overline{M}^{\com^*}_{g_1,g_2}]
\Big)|_{t=0} =
 \iota_*\Big( \sum_{a+b=2g+r-1}
(-1)^g (-1)^a 2g_2 \ \psi_{\star_1}^a \psi_{\star_2}^b
\scap[\Delta_{1,\emptyset}(g_1,g_2)] \Big)
$$
for $g_1,g_2 >0$.
Because of the $2g_2 t$ factor, the degenerate case $(g_1,g_2)=(g,0)$
contributes 0. However,
$$\epsilon_*\Big(
\frac{I_{g,r}}
{e(\text{Norm}_{0,g}^{vir})} \scap [\overline{M}^{\com^*}_{0,g}]
\Big)|_{t=0} =
(-1)^g (-1) 2g \ \psi_{1}^{2g+r}.$$
By the vanishing \eqref{vvp}, we conclude
$$(-1)^g (-1) 2g \ \psi_{1}^{2g+r} +
\sum_{g_1+g_2=g, \ g_i >0} \  \
\sum_{a+b=2g+r-1} \ \iota_*\Big(
(-1)^g (-1)^a 2g_2 \ \psi_{\star_1}^a \psi_{\star_2}^b
\scap[\Delta_{1,\emptyset}(g_1,g_2)]\Big)= 0$$
which is equivalent to Theorem \ref{bbt}. \qed

\subsection{Variations}
\label{nnvv}
Let $g\geq 0$ and $n_1,n_2 \geq 2$.
Consider the moduli
space $\overline{M}_{g,n_1+n_2}$.
Let $N_1$ and $N_2$ denote the markings sets
$$N_1=\{ 1, \ldots, n_1\}, \ \ N_2=\{n_1+1, \ldots, n_1+n_2\}.$$
For $g_1,g_2 \geq 0$, let
$$\iota: \DD_{N_1,N_2}[g_1,g_2]\rightarrow \overline{M}_{g,n_1+n_2}$$
denote the boundary divisor parametrizing
reducible curves $$C=C_1\cup C_2$$
with markings $N_i$ on $C_i$
satisfying $g(C_i)=g_i$ and $\ C_1 \cap C_2=p_\star$.
Let
$$\psi_{\star_1}, \psi_{\star_2} \in A^1\big(\DD_{N_1,N_2}(g_1,g_2)\big)$$
denote the cotangent line classes of $p_\star$ along
$C_1$ and $C_2$ as before.

\begin{prop} \label{bbbtt} For $g\geq 0$ and $n_1,n_2\geq 2$ and $r\geq 0$,
$$
\sum_{g_1+g_2=g,\ g_i\geq0}
\ \ \sum_{a+b=2g+n_1+n_2-3+r}
(-1)^a   \iota_*\big(\psi_{\star_1}^a \psi_{\star_2}^b \scap
[\DD_{N_1,N_2}(g_1,g_2)]\big)=0$$
in $A^{2g+n_1+n_2-2+r}(\overline{M}_{g,n_1+n_2})$.
\end{prop}

\begin{proof}
Consider the moduli space $\overline{M}_{g,n_1+n_2}(\proj^1,1)$
with the $\com^*$-action specified in Section \ref{axx}.
Let $$J_{g,r} =
\ev_1^*([\infty]^{1+r})\scup
\prod_{i\in N_1} \ev_i^*([\infty]) \scup
\prod_{i\in N_2} \ev_i^*([0]) \scup c_g(B) \in A^{g+n_1+n_2+r+1}
\big(\overline{M}_{g,n_1+n_2}(\proj^1,1)\big).$$
Since the non-equivariant limit of $[\infty]^2$ is 0, the non-equivariant
limit of $J_{g,r}$ is also 0.
Let
$$\epsilon: \overline{M}_{g,n_1+n_2}(\proj^1,1)
\rightarrow \overline{M}_{g,n_1+n_2}$$
be the forgetful map.
After push-forward,
\begin{equation}\label{iiirg}
\epsilon_*\big(
J_{g,r} \scap [\overline{M}_{g,n_1+n_2}(\proj^1,1)]^{vir}\big) \in
A_{\com^*}^{2g+n_1+n_2-2+r}(\overline{M}_{g,n_1+n_2}).
\end{equation}
Setting the non-equivariant limit to 0,
\begin{equation}\label{vvpc}
\epsilon_*\big(J_{g,r} \scap
[\overline{M}_{g,n_1+n_2}(\proj^1,1)]^{vir}\big)|_{t=0} = 0,
\end{equation}
yields an equation in $R^{2g+n_1+n_2-2+r}(\overline{M}_{g,n_1+n_2})$.
Evaluating the virtual localization formula as in the proof
of Theorem \ref{bbt} precisely yields Proposition \ref{bbbtt}.
\end{proof}

Since $n_1,n_2\geq 2$ in the hypothesis of Proposition \ref{bbbtt},
there are no degenerate cases. There is no difficulty in
handling the degenerate cases. We single out the following
result with the same proof{\footnote{
The proofs of Theorem 1 and Proposition 1 are almost
identical. In fact,
Theorem 1 can 
be derived from Proposition 1 using string and dilaton
equations.}}
 as Proposition \ref{bbbtt}.

\begin{prop}  \label{fqq}
For $g\geq 1$ and $r\geq 0$,
$$
-\psi_1^{2g+r}+ (-1)^r \psi_2^{2g+r} +
\sum_{g_1+g_2=g,\ g_i > 0}
\ \ \sum_{a+b=2g-1+r}
(-1)^a   \iota_*\big(\psi_{\star_1}^a \psi_{\star_2}^b \scap
[\DD_{1,2}(g_1,g_2)]\big)=0$$
in $A^{2g+r}(\overline{M}_{g,2})$.
\end{prop}

Proposition \ref{fqq} corresponds simply to the $n_1=n_2=1$
case of Proposition \ref{bbbtt}. The first two terms are
the degenerate contributions.

\subsection{Proof of Theorem \ref{vyt}}
We start by pushing forward the relation of Proposition
\ref{fqq} in genus $g+1$ to $\overline{M}_{g+1}$ for odd $r$,
$$
-2\kappa_{2(g+1)+r-2} -
\sum_{g_1+g_2=g+1,\ g_i > 0}
\ \ \sum_{a+b=2(g+1)-3+r}
(-1)^a   \iota_*\big(\psi_{\star_1}^a \psi_{\star_2}^b \scap
[\DD_{\emptyset,\emptyset}(g_1,g_2)]\big)=0,
$$
using the definition of the $\kappa$ classes and the
string equation.
Equivalently,
\begin{equation} \label{vpe}
\kappa_{2g+r} +\frac{1}{2}
\sum_{g_1+g_2=g+1,\ g_i > 0}
\ \ \sum_{a+b=2g-1+r}
(-1)^a   \iota_*\big(\psi_{\star_1}^a \psi_{\star_2}^b \scap
[\DD_{\emptyset,\emptyset}(g_1,g_2)]\big)=0 \in
A^{2g+r}(\overline{M}_{g+1})
\end{equation}
for odd $r$.

The Chern characters of the Hodge bundle
$\text{ch}_{2l-1}(\mathbb{E}_{g+1})$
on $\overline{M}_{g+1}$ vanish for $l> g+1$, see \cite{fp1}.
Hence, by Mumford's GRR calculation,
\begin{multline*}
\text{ch}_{2g+r}(\mathbb{E}_{g+1})
\left(\frac{B_{2g+r+1}}{(2g+r+1)!} \right)^{-1} = \\
\kappa_{2g+r} +\frac{1}{2} \iota_*(\xi_{g,r-1}) +    \frac{1}{2}
\sum_{g_1+g_2=g+1,\ g_i > 0}
\ \ \sum_{a+b=2g-1+r}
(-1)^a   \iota_*\big(\psi_{\star_1}^a \psi_{\star_2}^b \scap
[\DD_{\emptyset,\emptyset}(g_1,g_2)]\big)=0
\end{multline*}
for $r\geq 3$ odd.
Using the vanishing \eqref{vpe}, we conclude
$$\iota_*(\xi_{g,r-1}) = 0 \in A^{2g+r}(\overline{M}_{g+1})$$
for $r\geq 3$ odd, which are the only nontrivial cases of
Theorem \ref{vyt}.
\qed

\label{xxzz}

\section{Gromov-Witten equations}
\label{xl}


\subsection{Liu-Xu conjecture}
Let $X$ be a nonsingular projective variety.
We  prove here the
following result constraining the Gromov-Witten theory
of $X$ conjectured by K. Liu and H. Xu in
\cite{kliu}.

\begin{theorem} \label{conj:C}
Let $g\geq 0$ and
$x_{i}, y_{j} \in H^{*}(X, \mathbb{C})$. For all
 $p_{i}, q_{j}, r, s\geq 0$ and $m \geq 2g-3+r+s$,
\[ \sum_{k \in \mathbb{Z}} \
     \sum_{g_1+g_2=g, \ g_i\geq 0} \
(-1)^{k} \gwiih{g_1}{\graval{k} \prod_{i=1}^{r} \tau_{p_{i}}(x_{i})}
    \gwiih{g_2}{\gravual{m-k} \prod_{j=1}^{s} \tau_{q_{j}}(y_{j}) } = 0. \]
\end{theorem}

Here, $k$ is allowed to be an arbitrary integer.
To interpret  Theorem \ref{conj:C} correctly,
the following convention is used{\footnote{$\gamma_{1}$
is the identity of the cohomology ring of $X$.}}:
\begin{equation} \label{eqn:negdesc}
 \gwih{0,0}{\tau_{-2}(\gamma_{1})} = 1 \hspace{20pt} {\rm and} \hspace{20pt}
\gwih{0,0}{\grava{m} \gravb{-1-m}} = (-1)^{{\rm max}(m, -1-m)} \eta_{\alpha \beta}
\end{equation}
for $m \in \mathbb{Z}$. All other negative descendents vanish.
The sum over $\ell$ in Theorem \ref{conj:C} is implicit.

Since the genus 0 case of Theorem \ref{conj:C}
has been proved{\footnote{The $m \geq 3g-3+r+s$ case
is also proved in \cite{xliu2}, but the result will not be used here.}}
in \cite{xliu2}, we will
only consider the case $g \geq 1$.
By Theorem 0.2 of \cite{xliu2},
 Theorem~\ref{gwwq} follows from the $r=s=0$  case
of  Theorem \ref{conj:C}.

\subsection{Conventions}
We will not use convention \eqref{eqn:negdesc}.
Instead, we set $\grava{n}=0$ for $n<0$ and
separate the negative terms in the summation of  Theorem \ref{conj:C}.

The {\it big phase space} is the infinite dimensional vector space
with coordinate $t=(t_{n}^{\alpha})$. It can be interpreted as an
infinite product of the cohomology space $H^{*}(X,\com)$. The
Gromov-Witten potential $F_{g}^{X}$ is a function on the big phase
space. We will interpret the symbol $\grava{n}$ as the coordinate
vector field $\frac{\partial}{\partial t_{n}^{\alpha}}$. Moreover,
we also extend the meaning of $\gwiig{\vw_{1} \, \cdots \, \vw_{k}}$
from partial derivatives of $F_{g}^{X}$ to covariant derivatives of
$F_{g}^{X}$ with respect to arbitrary vector fields $\vw_{1},
\ldots, \vw_{k}$ on the big phase space. Here,
 the covariant differentiation
is with respect to the trivial connection $\nabla$ for which
the  coordinate vector fields
$\grava{n}$ are parallel.
More precisely, if
$\vw_{i} = \sum_{n, \alpha} f_{n, \alpha}^{i} \grava{n}$
where $f_{n, \alpha}^{i}$ are functions of $t=(t_{m}^{\beta})$, then
we define
\[ \gwiig{\vw_{1} \, \cdots \, \vw_{k}} =
\nabla^{k}_{\vw_{1}, \cdots, \vw_{k}} F_{g}^{X}
    = \sum_{\begin{array}{c} n_{1}, \cdots, n_{k} \\ \alpha_{1}, \cdots, \alpha_{k} \end{array}}
            \left( \prod_{i=1}^{k} f_{n_{i}, \alpha_{i}}^{i} \right)
            \gwiig{\tau_{n_{1}}(\alpha_{1}) \, \cdots \, \tau_{n_{k}}(\alpha_{k})}. \]

For a vector field of type $\grava{n}$, the integer
 $n$ is called the {\it level of the descendent}.
A vector field is {\it primary} if the level of the descendent is 0.
The total level of descendents for a set of vector fields
is defined to be the sum
of the levels of descendents for all vector fields in the set.
For convenience, we define the operators
$\tau_{+}$ and $\tau_{-}$ on the space of vector fields
by the following formulas:
\[ \tau_{\pm}(\vw) = \sum_{n, \alpha} f_{n, \alpha} \grava{n \pm 1}
\hspace{20pt} {\rm if} \hspace{20pt} \vw= \sum_{n, \alpha} f_{n, \alpha} \grava{n}.\]
Moreover, we define $\tau_{k}(\vw) = \tau_{+}^{k}(\vw)$
for any vector field $\vw$.

\subsection{Lower cases}
We first prove a result about relations among different cases of
Theorem \ref{conj:C}.
\begin{prop} \label{thm:reduceAB}
Let $g\geq 0$ be fixed.
If Theorem \ref{conj:C} holds for $r=\hat{r}$ and $s=\hat{s}$,
then  Theorem \ref{conj:C}
holds for all $r \leq \hat{r}$ and
$s \leq \hat{s}$.
\end{prop}

\noindent{\bf Proof}: We first rewrite Theorem \ref{conj:C} without
using the special convention \eqref{eqn:negdesc}.  Define
\[ \tilde{t}_{n}^{\alpha} = t_{n}^{\alpha} -\delta_{\alpha, 1}
\delta_{n, 1} .\] Let $\vw_{i}$, $\vv_{j}$ be arbitrary coordinate
vector
 fields on the big phase space of the form $\grava{n}$.
For $r, s, g, m \geq 0$, define
\begin{multline} \label{eqn:Psi}
\Psi_{r,s, g, m} (\vw_{1}, \cdots , \vw_{r} \mid \vv_{1}, \cdots, \vv_{s})= \\
\sum_{k=0}^{m} \ \ \sum_{g_1+g_2=g, \ g_i\geq 0} \
           (-1)^k \gwiih{g_1}{\grava{k} \, \vw_{1} \, \cdots \, \vw_{r}}
                \gwiih{g_2}{\gravua{m-k} \, \vv_{1} \, \cdots \, \vv_{s}}
                \\
 - \delta_{r,0}  \sum_{n, \alpha} \tilde{t}_{n}^{\alpha}
                \gwiig{\grava{n+m+1} \, \vv_{1} \, \cdots \, \vv_{s}}
- \delta_{r,1} \gwiig{\tau_{m+1}(\vw_{1}) \, \vv_{1} \, \cdots \, \vv_{s})}
     \\
+ \delta_{s,0} (-1)^{m+1} \sum_{n, \alpha} \tilde{t}_{n}^{\alpha}
                \gwiig{\grava{n+m+1} \, \vw_{1} \, \cdots \, \vw_{r}}
+ \delta_{s,1} (-1)^{m+1}
                \gwiig{\vw_{1} \, \cdots \, \vw_{r} \, \tau_{m+1}(\vv_{1})}.    \end{multline}
The function satisfies the symmetry
\begin{equation} \label{eqn:psir->s}
\Psi_{r,s, g, m}(\vw_{1} \cdots  \vw_{r} \mid \vv_{1} \cdots \vv_{s})
= (-1)^{m} \Psi_{s,r, g, m}(\vv_{1} \cdots  \vv_{s} 
\mid \vw_{1} \cdots \vw_{r}).
\end{equation}
Moreover, $\Psi_{0, 0, g, m}$ is identically equal to 0 if $m$ is odd.

Theorem \ref{conj:C} can be restated as
\begin{equation} \label{eqn:conjC}
\Psi_{r, s, g, m}(\vw_{1} \cdots  \vw_{r} \mid \vv_{1} \cdots \vv_{s}) = 0
\end{equation}
if $m \geq 2g+r+s-3$, see \cite{kliu2}.


Suppose for fixed integers $r>0$ and $s \geq 0$,
equation~\eqref{eqn:conjC} holds  for all integers $m \geq 2g+r+s-3$.
Then, we must prove that
equation~\eqref{eqn:conjC}
holds if $r$ is replaced by $r-1$ for all
$m \geq 2g+r+s-4$.
By an inverse induction on $r$, if Theorem \ref{conj:C} holds for $r=\hat{r}$
and $s=\hat{s}$, then Theorem \ref{conj:C}  holds for $r \leq \hat{r}$ and $s=\hat{s}$.
By equation \eqref{eqn:psir->s}, we can switch the role of $r$ and $s$.
Hence, the Proposition will be proved.

Consider the {\em string vector field},
\[ \vs= - \sum_{n, \alpha} \tilde{t}_{n}^{\alpha} \grava{n-1} .\]
The {\it string equation} for Gromov-Witten invariants can be written as
\[ \gwiig{\vs} = \frac{1}{2} \delta_{g, 0} \eta_{\alpha \beta}
t_{0}^{\alpha} t_{0}^{\beta} \]
where $\eta_{\alpha \beta}= \int_{X} \ga \cup \gb$ is the usual pairing.
Taking derivatives of the string equation, we obtain
\begin{equation} \label{eqn:DerString}
 \gwiig{ \vs \, \vw_{1} \, \cdots \, \vw_{k} }
    = \sum_{i=1}^{k} \gwiig{ \vw_{1} \, \cdots \,
            \left\{ \tau_{-}(\vw_{i}) \right\}
            \, \cdots \, \vw_{k} }
            + \delta_{g, 0} \nabla^{k}_{\vw_{1}, \cdots, \vw_{k}}
                \left( \frac{1}{2} \eta_{\alpha \beta}
            t_{0}^{\alpha} t_{0}^{\beta} \right).
\end{equation}
Note that
\[  \nabla^{k}_{\vw_{1}, \cdots, \vw_{k}}
                \left( \frac{1}{2} \eta_{\alpha \beta}
            t_{0}^{\alpha} t_{0}^{\beta} \right) = 0 \]
if $k>2$ or if at least one of the vector fields $\vw_{1}, \cdots, \vw_{k}$
has a positive descendent level.

Since equation~\eqref{eqn:conjC} is linear with respect to
 each $\vw_{i}$ and $\vv_{j}$, we can replace
them by any vector fields on the big phase space. Assume $r>0$.
We consider what happens if $\vw_{r} = \vs$ in
$\Psi_{r,s, g, m} (\vw_{1}, \cdots , \vw_{r} \mid \vv_{1}, \cdots, \vv_{s})$.

\begin{lemma} \label{lem:Sreduce}
For $r >0$,
\begin{multline}
\Psi_{r,s, g, m} (\vw_{1} \cdots  \vw_{r-1} \vs \mid \vv_{1} \cdots \vv_{s})
        \label{eqn:reducepsia3}  =\\
  - \Psi_{r-1,s, g, m-1} (\vw_{1} \cdots  \vw_{r-1} \mid \vv_{1} 
\cdots \vv_{s}) \\
 + \sum_{i=1}^{r-1} \Psi_{r-1,s, g, m} (\vw_{1} \cdots  \tau_{-}(\vw_{i})
                        \cdots \vw_{r-1} \mid \vv_{1} \cdots \vv_{s}) .
\end{multline}
for all vector fields $\vw_{i}$ and $\vv_{j}$.
\end{lemma}

Assuming the validity of Lemma \ref{lem:Sreduce},
we can prove the Proposition by induction. Indeed,
assume
\[ \Psi_{r,s, g, m} (\vw_{1} \cdots  \vw_{r-1} \vs \mid \vv_{1} 
\cdots \vv_{s}) = 0, \]
for all vector fields $\vw_{i}$, $\vv_{j}$, and all integers $m \geq 2g-3+r+s$.
By linearity, we may assume that all vector fields $\vw_{i}$ are coordinate vector fields
of type $\grava{n}$.
Note that $\tau_{-}(\vw_{i}) = 0$ if $\vw_{i}$ is a primary vector field.
Hence equation~\eqref{eqn:reducepsia3} implies
\begin{equation} \label{eqn:psia-1a3}
 \Psi_{r-1,s, g, m-1} (\vw_{1} \cdots \vw_{r-1} \mid \vv_{1}
 \cdots \vv_{s})=0
\end{equation}
for all integers $m \geq 2g-3+r+s$
if $\vw_{1}, \cdots, \vw_{r-1}$ are all primary vector fields.

Since the total level of descendents for vector fields in the second term on the
right hand side of equation~\eqref{eqn:reducepsia3} is strictly less than
that in the first term, an induction on the total level of descendents for
$\vw_{1}, \cdots, \vw_{r-1}$ shows that equation~\eqref{eqn:psia-1a3} also hold for
all (not necessarily primary) vector fields $\vw_{1}, \cdots, \vw_{r-1}$.
Hence, if Theorem \ref{conj:C} holds for $r>0$ and $s \geq 0$, then Theorem \ref{conj:C}
holds if $r$ is replaced by $r-1$.  The Proposition thus
 follows from Lemma~\ref{lem:Sreduce}.\qed

\subsection{Proof of  Lemma~\ref{lem:Sreduce}}
 Using equation~\eqref{eqn:DerString},
the result is  straightforward for $r>2$.
The cases $r \leq 2$ are more subtle
because of  the last term in equation~\eqref{eqn:DerString}.

We consider the case $r=2$ first. If $\vw$ is a primary vector field, then
\[ \nabla^{2}_{\vw, \grava{k}}
                \left( \frac{1}{2} \eta_{\beta \mu}
            t_{0}^{\beta} t_{0}^{\mu} \right) \gwiig{\gravua{m-k} \, \vv_{1} \, \cdots \, \vv_{s}}
            \,\,=\,\, \delta_{k, 0} \gwiig{\tau_{m}(\vw) \, \vv_{1} \, \cdots \, \vv_{s}}. \]
This will produce the extra term in $\Psi_{1,s, g, m-1} (\vw  \mid \vv_{1}, \cdots, \vv_{s})$.
Therefore by equation~\eqref{eqn:DerString},
\begin{eqnarray}
 \Psi_{2,s, g, m} (\vw  \vs \mid \vv_{1} \cdots \vv_{s})
&=&  - \Psi_{1,s, g, m-1} (\vw  \mid \vv_{1} \cdots \vv_{s})
    \nonumber
\end{eqnarray}
when $\vw$ is a primary vector field.

If $\vw$ has a positive descendent level, then
\[ \nabla^{2}_{\vw, \grava{k}}
                \left( \frac{1}{2} \eta_{\beta \mu}
            t_{0}^{\beta} t_{0}^{\mu} \right) = 0 \]
for all $k\geq 0$. Equation~\eqref{eqn:DerString} again implies
\begin{multline}
 \Psi_{2,s, g, m} (\vw \vs \mid \vv_{1} \cdots \vv_{s})
        \label{eqn:psia-1a2} = \\
  - \Psi_{1,s, g, m-1} (\vw \mid \vv_{1} \cdots \vv_{s})
 +  \Psi_{1,s, g, m} (\tau_{-}(\vw) \mid \vv_{1} \cdots \vv_{s}) .
\end{multline}
When passing from $\Psi_{2,s, g, m}$ to $\Psi_{1,s, g, m}$, an extra term
will emerge. The summations which we obtain from applying
equation~\eqref{eqn:DerString} to
$ \Psi_{2,s, g, m} (\vw, \vs \mid \vv_{1}, \cdots, \vv_{s}) $
have some missing terms when compared to the definition of
$\Psi_{1,s, g, m-1} (\vw \mid \vv_{1} \cdots \vv_{s})$ and
$\Psi_{1,s, g, m} (\tau_{-}(\vw) \mid \vv_{1} \cdots \vv_{s})$.
The missing term for $\Psi_{1,s, g, m-1} (\vw \mid \vv_{1} \cdots \vv_{s})$
is $- \gwiig{\tau_{m}(\vw) \, \vv_{1} \cdots \vv_{s}}$ while the missing term
for $\Psi_{1,s, g, m} (\tau_{-}(\vw) \mid \vv_{1}\, \cdots\, \vv_{s})$
is $- \gwiig{\tau_{m+1}(\tau_{-}(\vw)) \, \vv_{1}\, \cdots\, \vv_{s}}$.
The missing terms cancel
in \eqref{eqn:psia-1a2}
when $\vw$ has a positive descendent level.

Since we have checked that equation~\eqref{eqn:psia-1a2}
holds for all primary and descendent vector fields $\vw$,
Lemma~\ref{lem:Sreduce} is true for $r=2$.

Consider next  the case $r=1$ and $s>0$. We have
\begin{multline*}
\Psi_{1,s, g, m}(\vs \mid \vv_{1} \cdots \vv_{s}) = \\
 - \gwiig{\tau_{m+1}(\vs) \, \vv_{1} \, \cdots \, \vv_{s}}
    + \delta_{s, 1} (-1)^{m+1} \gwiig{\vs \,\, \tau_{m+1}(\vv_{1})} \\
    + \sum_{k=0}^{m}\ \ \sum_{g_1+g_2=g,\ g_i\geq 0}\ (-1)^{k} 
\gwiih{g_1}{\tau_{k}(\ga) \, \vs}
        \gwiih{g_2}{\gravua{m-k} \, \vv_{1} \, \cdots \, \vv_{s}}.
\end{multline*}
In the definition of $\vs$, $\tilde{t}_{0}^\alpha$ is not included since
$\tau_{-1}(\ga)=0$. Hence
\begin{equation} \label{eqn:tau+S}
 \gwiig{\tau_{m+1}(\vs) \, \vv_{1} \, \cdots \, \vv_{s}}
 = - \sum_{n=1}^{\infty} \sum_{\alpha} \tilde{t}_{n}^{\alpha}
    \gwiig{\grava{n+m} \, \vv_{1} \, \cdots \, \vv_{s}}.
\end{equation}
By equation~\eqref{eqn:DerString},
$\gwiig{\vs \,\, \tau_{m+1}(\vv_{1})} = \gwiig{\tau_{m}(\vv_{1})}$ and
\[ \gwiih{g_1}{\tau_{k}(\ga) \, \vs} = \gwiih{g_1}{\tau_{k-1}(\ga)} +
            \delta_{g_1, 0} \delta_{k, 0} \eta_{\alpha \beta} t_{0}^{\beta}.\]
The effect of the second term on the right hand side of this equation
is just to compensate for the missing case $n=0$ in the summation
for $n$ in equation~\eqref{eqn:tau+S} when computing
$$\Psi_{1,s, g, m}(\vs \mid \vv_{1}, \cdots, \vv_{s}).$$
Therefore we have
\begin{eqnarray*}
  \Psi_{1,s, g, m}(\vs \mid \vv_{1} \cdots \vv_{s})
&=& - \Psi_{0,s, g, m-1}(\vv_{1} \cdots \vv_{s}).
\end{eqnarray*}
Hence, Lemma~\ref{lem:Sreduce} is true for $r=1$ and $s>0$.

Now only the case $r=1$ and $s=0$ is left. By definition,
\begin{eqnarray*}
\Psi_{1,0, g, m}(\vs)
&=& - \gwiig{\tau_{m+1}(\vs)}
    + (-1)^{m+1} \sum_{n, \alpha} \tilde{t}_{n}^{\alpha}
        \gwiig{\grava{n+m+1} \,\, \vs} \\
&&    +  \sum_{k=0}^{m}\ \ 
\sum_{g_1+g_2=g,\ g_i\geq 0}\ (-1)^{k} \gwiih{g_1}{\tau_{k}(\ga) \, \vs}
        \gwiih{g_2}{\gravua{m-k} }.
\end{eqnarray*}
By equation~\eqref{eqn:DerString}, we have
\begin{eqnarray*}
\Psi_{1,0, g, m}(\vs)
&=&  \left\{1+(-1)^{m+1} \right\} \sum_{n, \alpha} \tilde{t}_{n}^{\alpha}
        \gwiig{\grava{n+m} } \\
&&    - \sum_{k=0}^{m-1}\ \ \sum_{g_1+g_2=g,\ g_i\geq 0} \
 (-1)^{k} \gwiih{g_1}{\tau_{k}(\ga) }
        \gwiih{g_2}{\gravua{m-1-k} } \\
&=& -\Psi_{0,0, g, m-1}.
\end{eqnarray*}
The proof for Lemma~\ref{lem:Sreduce} is complete. \qed

\subsection{Proof of Theorem \ref{conj:C}}
Relations in  $R^{*}(\barr{M}_{g, n})$
can be translated into universal equations for Gromov-Witten invariants
by the splitting axiom and cotangent line comparison equations.
Define the operator $T$ on the space of vector fields by
\[ T(\vw) = \tau_{+}(\vw) - \gwii{\vw \, \gua} \ga \]
for any vector field $\vw$. Properties of $T$ have been studied in \cite{xliu1}.
The operator is very useful for the translation into universal equations.
In the process, each marked point corresponds to a vector field, and the cotangent line
class corresponds to the operator $T$. Each node is translated into a pair of
primary vector fields $\gamma_\ell$ and $\gamma^\ell$. In particular, the relation
of Proposition~\ref{bbbtt} is translated into the following universal
equation
\begin{equation} \label{eqn:rs2T}
  \sum_{k=0}^{m}\ \ \sum_{g_1+g_2=g, \ g_i\geq 0} \  (-1)^{k}
    \gwiih{g_1}{\vw_{1} \, \cdots \, \vw_{n_{1}} \, T^{k}(\gamma_\ell)}
    \gwiih{g_2}{T^{m-k}(\gamma^\ell) \, \vv_{1} \, \cdots \, \vv_{n_{2}}} = 0
\end{equation}
for all vector fields $\vw_{i}$ and $\vv_{j}$ if $n_{1}, n_{2} \geq 2$ and $m \geq 2g+n_{1}+n_{2}-3$.

Let $P$ and $Q$ be two arbitrary contravariant tensors on the big phase space.
The following formula was proved in \cite[Proposition 3.2]{xliu2}:
\[ \sum_{k=0}^{m} (-1)^{k} P(T^{k}(\gamma_\ell)) \, \, Q(T^{m-k}(\gamma^\ell))
    = \sum_{k=0}^{m} (-1)^{k} P(\graval{k}) \,\,  Q(\gravual{m-k}) \]
for $m \geq 0$.
In particular, if we take 
$P(\vu)= \gwiih{g_1}{\vw_{1} \, \cdots \, \vw_{n_{1}} \, \vu}$
and $Q(\vu)= \gwiih{g_2}{\vu \, \, \vv_{1} \, \cdots \, \vv_{n_{2}}}$, then the left hand side of
equation~\eqref{eqn:rs2T} is equal to
\begin{multline*} 
\sum_{k=0}^{m} \ \ \sum_{g_1+g_2=g, \ g_i\geq 0} \
 (-1)^{k}
    \gwiih{g_1}{\vw_{1} \, \cdots \, \vw_{n_{1}} \, \tau_{k}(\gamma_\ell)}
    \gwiih{g_2}{\tau_{m-k}(\gamma^\ell) \, \vv_{1} \, \cdots \, \vv_{n_{2}}} =\\
 \Psi_{n_{1}, n_{2}, g, m}(\vw_{1} \cdots \, \vw_{n_{1}} \mid
                 \vv_{1} \cdots \vv_{n_{2}}).
\end{multline*}
Therefore equation~\eqref{eqn:rs2T} implies that Theorem ~\ref{conj:C} is true
for $r=n_{1} \geq 2$
and $s = n_{2} \geq 2$. By Proposition \ref{thm:reduceAB}, all other cases of Theorem
\ref{conj:C}
follow. \qed

\end{document}